# Sharp failure rates for the bootstrap particle filter in high dimensions

## Peter Bickel[1], Bo Li[2] and Thomas Bengtsson[3]

*University of California-Berkeley, Tsinghua University and Bell Labs*

**Abstract:** We prove that the maximum of the sample importance weights in a high-dimensional Gaussian particle filter converges to unity unless the ensemble size grows exponentially in the system dimension. Our work is motivated by and parallels the derivations of Bengtsson, Bickel and Li (2007); however, we weaken their assumptions on the eigenvalues of the covariance matrix of the prior distribution and establish rigorously their strong conjecture on when weight collapse occurs. Specifically, we remove the assumption that the nonzero eigenvalues are bounded away from zero, which, although the dimension of the involved vectors grow to infinity, essentially permits the effective system dimension to be bounded. Moreover, with some restrictions on the rate of growth of the maximum eigenvalue, we relax their assumption that the eigenvalues are bounded from above, allowing the system to be dominated by a single mode.

## Contents

## 1. Introduction

Bayesian filtering methods are a commonly employed tool for combing physical models and data. The filters treat the unknown system state as a random variable and resolve its probability density conditional on the data (and the system dynamics) through Monte Carlo sampling techniques. When applied sequentially in time, these methods are commonly referred to as particle filters ([8], [10]). For a diverse collection of applications and an excellent introduction to the field in general, see the edited volume by Doucet [6]. The particle filter method relies heavily on a likelihood based reweighting mechanism of the involved sample draws. This reweighting scheme produces the so called importance weights, and these weights are the primary focus of our work. Specifically, in a Gaussian filter context, we examine the behavior of the importance weights as a function of the system dimension and of sample size.

---

[1]367 Evans Hall, Department of Statistics, University of California, Berkeley, 94710-3860, CA, USA, e-mail: bickel@stat.berkeley.edu

[2]440 Weilun Hall, School of Economics and Management, Tsinghua University, Beijing, 100084, China, e-mail: libo@sem.tsinghua.edu

[3]Bell Labs, 600 Mountain Avenue, Murray Hill, 07904 NJ, USA, e-mail: tocke@cgd.ucar.edu
*AMS 2000 subject classifications:* 93E11, 62L12, 86A22, 60G50, 86A32, 86A10.
*Keywords and phrases:* Bayesian filter, curse of dimensionality, ensemble forecast, ensemble methods, importance sampling, large deviations, Monte Carlo, numerical weather prediction, sample size requirements, state-space model.

The popularity of the particle filter is no doubt due to the flexibility of the model framework to handle both non-linear and non-gaussian structures. However, in spite of its generality, the method is not without flaws: the particle filter is known to require large Monte Carlo ensembles and frequent resampling to estimate the desired densities (cf., [9]). This drawback is particularly prevalent in higher dimensional systems where the filter becomes unstable and quickly collapses onto a single point mass. In recent work, for a single Bayes update step in a Gaussian setting, Bengtsson, Bickel, and Li [3] give a derivation of the weight collapse as a function of the system dimension and of sample size. To shed further light on the weight collapse, this paper establishes conjectures (given in [3]) which make their arguments fully rigorous. Just as significantly, we exhibit that collapse is a function of the *effective dimension* (defined in Section 3), rather than the absolute dimension. As in [3], our analysis is given in the context of a stylized Gaussian example, but we conjecture (and simulations show) that our results are informative for situations that depend on similarly defined reweighting schemes. The results imply that to avoid collapse, the sample size must grow super-exponentially in the effective dimension. We do not investigate refinements of particle filters methods, such as simulated tempering [4], although our discussion in Section 2.1 suggests that their approach is not a solution to avoid collapse in truly high-dimensional settings.

Our work is outlined as follows. The next section describes the particle filter, provides notation, and describes the use of the ensemble method for approximating posterior densities. The main developments are then presented in Section 3, where we give several results establishing the conditions under which the maximum sample weight in a Gaussian particle filter converges to unity. All technical results are proved in the Appendix. (We note that some material in Section 2.1 and Section 3 is given in [3], but is reproduced here for completeness.)

## 2. Model setting

### 2.1. *The particle filter*

Let $X_t$ represent the unknown system state at time $t$, $Y_t$ be a noisy data measurement of $X_t$, and let $\mathbf{Y}^t$ represent all data up to and including time $t$. Based on the data $\mathbf{Y}^t$ and (some) knowledge of the time-evolution of the system state from $X_{t-1}$ to $X_t$, we seek the posterior distribution $p(X_t|\mathbf{Y}^t)$. We assume we have available a random sample $\{X_{t,i}^f\}$ of size $n$ from the prior distribution $p(X_t|\mathbf{Y}^{t-1})$. Associated with the prior sample is a set of weights $\{w_i^f\}$. We assume further that the likelihood density $p(Y_t|X_t)$ is computable for arbitrary $X_t$.

The particle filter seeks to recursively in time estimate the probability distribution of the unknown state $X_t$. At each time $t$, the probability distribution is represented by the sample ensemble $\{X_{t,i}^f, w_i^f\}$, and the distribution can be propagated forward one time-step by evolving each $X_{t,i}^f$ using the system dynamics. Once new data $Y_t$ is available, Bayes theorem is used to adjust the weights based on how

"close" the associated sample points are to the data. The following schematic describes the particle filter:

$$p(X_t|\mathbf{Y}^{t-1}), Y_t \xrightarrow{\text{Bayes}} p(X_t|\mathbf{Y}^t) \xrightarrow{G(\cdot)} p(X_{t+1}|\mathbf{Y}^t), Y_{t+1} \xrightarrow{\text{Bayes}} p(X_{t+1}|\mathbf{Y}^{t+1}).$$

Here, at time $t$ (on the left), Bayes theorem combines $p(X_t|\mathbf{Y}^{t-1})$ and $Y_t$ to produce $p(X_t|\mathbf{Y}^t)$. The system dynamics, in the above represented by $G(\cdot)$ (middle), is used to propagate $p(X_t|\mathbf{Y}^t)$ one time step and this yields $p(X_{t+1}|\mathbf{Y}^t)$. Bayes theorem is then again employed to find the posterior $p(X_{t+1}|\mathbf{Y}^{t+1})$ (right).

In a particle filter, the above schematic is straightforwardly implemented (at least conceptually) using a random sample. We note first that the change-of-variables problem represented by the propagation of $p(X_t|\mathbf{Y}^t)$ can be solved by evaluating $G(\cdot)$ at each sample point. We will not discuss the implementation of the forecast step here; instead, our focus is on the Bayes update step. As mentioned, the particle filter implements the Bayes step by reweighting the prior sample according to the likelihood. We note in passing that the particle filter may be derived as a (sequential) importance sampler (e.g., [2]) where the proposal distribution is given by the prior and the target distribution is given by the posterior. In the schematic below, which describes a bootstrap-likelihood filter, the prior sample is "converted" to a posterior sample by resampling (with replacement) each member $X_{t,i}^f$ with probability proportional to $w_i^f \times p(Y_t|X_{t,i}^f)$, i.e.,

$$\overbrace{\{X_{t,1}^f, \ldots, X_{t,n}^f\}}^{prior\ ensemble},\ Y_t \xrightarrow{resample} \overbrace{\{X_{t,1}^u, \ldots, X_{t,n}^u\}}^{posterior\ ensemble}.$$

Although the particle filter has been successfully applied to a variety settings, it often produces highly varying importance weights. Remedies to stabilize the filter include resampling (renormalizing) the involved empirical measure at regular time intervals [8, 9] and marginalizing or restricting the sample space by conditioning on a larger information set [10, 11]. Another approach is given by simulated tempering [4], which makes use of the regularized likelihood $p(Y_t|X_{t,i}^f)^\alpha$, where $0 < \alpha < 1$. However, as can be seen from our derivations, e.g. Proposition 3.1, a fixed $\alpha$ does not alter the conclusion of collapse. Moreover, for each time point, to obtain samples from the target density, simulated tempering generates a sequence of ensembles from kernels $K_i(\cdot)$ ($i = 1, \ldots, I$) such that $K_I(\cdot)$ approaches the desired kernel $K(\cdot)$ associated with the posterior density. Unfortunately, for truly high dimensional systems, we conjecture that the number of intermediate sampling steps $I$ would be prohibitively large and render it practically unfeasible. Thus, such remedies do not fundamentally address performance when the filter is applied to *very large scale* systems. For example, as noted by ([1], [13]), when applied in high dimensions, the filter collapses to a point mass after a few (or even one!) observation cycles. In particular, as will be shown in Section 3, it is the normalized quantity $w_i = p(Y_t|X_{t,i}^f)/\sum_j p(Y_t|X_{t,j}^f)$ that behaves singularly.

The next section sets up the necessary notation and formalizes our problem.

## 2.2. Monte Carlo scheme

We formalize our problem as follows. Consider a set of $n$ sample points $\mathbf{X} = \{X_1, \ldots, X_n\}$, where $X_i \in \Re^d$ and both the sample size $n$ and system dimension $d$ are "large." (To lighten notation, we have dropped the time subscript and the forecast superscript.) We assume that the sample $\mathbf{X}$ is drawn randomly from the prior (or proposal) distribution $p(X)$. New data $Y$ is related to the state $X$ by the conditional density $p(Y|X)$. For concreteness, a functional relationship $Y = f(X) + \varepsilon$ is assumed, and $\varepsilon$ is taken to be independent of the state $X$. The goal is to estimate posterior expectations using the importance ratio, i.e., for some function $h(\cdot)$, we want to estimate

$$E(h(X)|Y) = \int h(X) \frac{p(Y|X)p(X)}{\int p(Y|X)p(X)\mathrm{d}X} \mathrm{d}X,$$

and use

$$\hat{E}(h(X)|Y) = \sum_{i=1}^{n} h(X_i) \frac{p(Y|X_i)}{\sum_{j=1}^{n} p(Y|X_j)}$$

as an estimator. Based on this formulation, the weights attached to each ensemble member

(1) $$w_i = \frac{p(Y|X_i)}{\sum_{j=1}^{n} p(Y|X_j)}$$

are the primary objects of our study. As mentioned, in large scale analyzes, the weights in (1) are highly variable and often produce estimates $\hat{E}(\cdot)$ which are collapsed onto a point mass with $max(w_i) \approx 1$. As illuminated in [3], this degeneracy is pervasive for high-dimensional systems, and appears to hold for a variety of prior and likelihood distributions.

We next consider the case when both the prior and the likelihood distributions are Gaussian.

## 3. Gaussian case

We assume a data model given by $Y = HX + \varepsilon$, where $Y$ is a $d \times 1$ vector, $H$ is a known $d \times q$ matrix, and $X$ is a $q \times 1$ vector. Both the proposal distribution and the error distribution are Gaussian with $p(X) = N(\mu_X, \Sigma_X)$ and $p(\varepsilon) = N(0, \Sigma_\varepsilon)$, and the noise $\varepsilon$ is taken independent of the state $X$. Since the data model can be pre-rotated by $\Sigma_\varepsilon^{-1/2}$, we set $\Sigma_\varepsilon = I_d$ without loss of generality (wlog). Moreover, since $EY = EHX$, we can replace $X_i$ by $(X_i - EX_i)$ and $Y$ by $(Y - EY)$ and leave $p(Y|X)$ unchanged. Hence, wlog we also set $\mu_X = 0$. Further, define, for conformable $A$ and $B$, the inner product $\langle A, B \rangle = A^T B$ (where the superscript $^T$ denotes matrix transpose), and let $\|A\|^2 = \langle A, A \rangle$.

With $p(Y|X) \sim N(HX, I_d)$, the weights in (1) can be expressed as

(2) $$w_i = \frac{\exp\big(-\|Y - HX_i\|^2/2\big)}{\sum_{j=1}^{n} \exp\big(-\|Y - HX_j\|^2/2\big)}.$$

To establish weight collapse for high-dimensional Gaussian $p(Y|X)$ and $p(X)$, we first write the exponent in (2) in terms of the singular values of $cov(HX)$.

Let $d' = rank(H)$. With $\lambda_1^2, \dots, \lambda_{d'}^2$ the singular values of $cov(HX)$, define the $d' \times d'$ matrix $D = diag(\lambda_1, \dots, \lambda_{d'})$. Then, with $Q$ an orthogonal matrix obtained by the singular value decomposition of $cov(HX)$, define the $d' \times 1$ vector $V$ such that

$$Q^T HX = DV.$$

Note that $V_i$ corresponding to $X_i$ is $N(0, I_{d'})$. Since $Q$ is orthogonal, we can write

$$(3) \qquad \|Y - HX_i\|^2 = \|Q^T Y - DV_i\|^2 = \sum_{j=1}^{d'} \lambda_j^2 W_{ij}^2 + \sum_{j=d'+1}^{d} \epsilon_{0j}^2,$$

where, conditional on $Y$, $[W_{i1}, \ldots, W_{id'}]^T$ is $N(\xi, I_{d'})$, and where $\epsilon_{0j}$ is the $j$th component of the observation noise vector $\varepsilon$. The mean vector $\xi = [\mu_1, \ldots, \mu_{d'}]^T$ is given by

$$(4) \qquad \xi = D^{-1} Q^T Y = V + D^{-1} \varepsilon',$$

where $V$ and $\varepsilon'$ are independent $N(0, I_{d'})$.

Now, for $i = 1, \ldots, n$, define

$$(5) \qquad S_i = \frac{\sum_{j=1}^{d'} \lambda_j^2 (W_{ij}^2 - (1 + \mu_j^2))}{\left(2 \sum_{j=1}^{d'} \lambda_j^4 (1 + 2\mu_j^2)\right)^{1/2}}.$$

Note that the second term in (3) is constant for every $X_i$, and will not influence the weight $w_i$.

By (2), we can express the maximum weight as

$$(6) \qquad w_{(n)} = \frac{1}{1 + T_{n,d'}},$$

where $T_{n,d'} = \sum_{\ell=2}^{n} e^{-\sigma_{d'}\sqrt{d'}(S_{(\ell)} - S_{(1)})}$ with $\sigma_{d'}^2 = \frac{2}{d'} \sum_{j=1}^{d'} \lambda_j^4 (1 + 2\mu_j^2)$. Thus, to prove weight collapse, we need to show convergence of the denominator in (6) to unity. We now state the following.

**Proposition 3.1.** *Let $S_i, i = 1, \ldots, n$, be independent random variables with cumulative distribution function (cdf) $G_d(\cdot)$ satisfying the conditions specified in Lemma A.1 and Lemma A.2 stated in the Appendix. Let $S_{(1)} \leq \cdots \leq S_{(n)}$ be the ordered sequence of $S_1, \ldots, S_n$, and define, for some $\sigma > 0$, $T_{n,d} = \sum_{\ell=2}^{n} e^{-\sigma\sqrt{d}(S_{(\ell)} - S_{(1)})}$. Then, as, $n, d \to \infty$, if $\frac{\log n \log d}{d} \to 0$, we have*

$$\sqrt{\frac{\sigma^2 d}{2 \log n}} E(T_{n,d}) \to 1.$$

A proof of the result is provided in the Appendix. For the Gaussian case considered here, an immediate implication of Proposition 3.1 is weight collapse. Specifically, with two additional assumptions, we may assert the following.

**Proposition 3.2.** *We assume, for the Gaussian case considered here,*

A1: *There is a positive constant $\delta$ such that $\frac{1}{\delta} \geq \lambda_1, \cdots, \lambda_{d'} > \delta$; and*

A2: $\tau_{d'}^2 = \frac{2}{d'} \sum_{j=1}^{d'} (3\lambda_j^4 + 2\lambda_j^2) \to \sigma^2 > 0$.

*Then, if $\frac{\log n \log d'}{d'} \to 0$, we have $w_{(n)} \xrightarrow{P} 1$.*

Proposition 3.2 follows by Lemma A.3 (Appendix) and Proposition 3.1.

The above result implies that, unless $n$ grows super-exponentially in $d'$, we have weight collapse. We note that Proposition 3.2 is a sharpening of the convergence rate as compared to that implied by Section 3.1 of [3]. The $\log d'$ term appears only because $\max |\mu_j| = O_p(\sqrt{\log d'})$, and we need to make our analysis conditional on the $\{\mu_j\}$.

The results in Proposition 3.2 suggest that large $d'$ leads to collapse. However, we argue now that what really matters is the *effective dimension* of $X$, defined as the sum of the singular values of $cov(HX)$. We shall assume that

B : $\lambda_1 \geq \lambda_2 \geq \cdots \geq \lambda_{d'} \geq \cdots$ are part of an infinite sequence.

Our arguments can be modified to the case where $\{\lambda_j : 1 \leq j \leq d'\}$ is a double array, but we eschew this complication.

There are two possibilities,

$$\text{(i)} \sum_{j=1}^{\infty} \lambda_j^2 < \infty, \ \text{ or (ii) } \sum_{j=1}^{\infty} \lambda_j^2 = \infty.$$

We claim that if (i) holds, there is no weight collapse. That is, if, say, $g : \mathcal{R} \mapsto \mathcal{R}$ is bounded and continuous,

$$(7) \qquad \sum_{i=1}^{n} w_i g(X_i^*) \xrightarrow{P} Eg(X|Y).$$

In the above, $X_i^*$ is drawn from the empirical measure $\sum_{j=1}^{n} w_i \delta(X_i)$, where $\delta(\cdot)$ represents the delta function, and where, as before, the $w_i$ represents the likelihood-defined weights.

To verify the convergence in (7), note that

$$w_i = U_i / \sum_{j=1}^{n} U_j,$$

where

$$(8) \qquad U_i = c_{d'}^{-1} \exp\{-\frac{1}{2} \sum_{j=1}^{d'} \left[\lambda_j^2(Z_{ij}^2 - 1) + 2\lambda_j^2 \mu_j Z_{ij}\right]\}.$$

In (8), the $Z_{ij}$'s are i.i.d. $N(0, 1)$ and

$$\begin{aligned}
c_{d'} &= E\left[\exp\{-\frac{1}{2} \sum_{j=1}^{d'} \left[\lambda_j^2(Z_{ij}^2 - 1) + 2\lambda_j^2 \mu_j Z_{ij}\right]\}\right] \\
&= \prod_{j=1}^{d'} \left[(1 + \lambda_j^2)^{-1/2} e^{\lambda_j^2/2} e^{\frac{\lambda_j^4 \mu_j^2}{2(1+\lambda_j^2)}}\right].
\end{aligned}$$

Now, since (i) implies that $\prod_{j=1}^{d'}(1 + \lambda_j^2)^{-1/2} e^{\lambda_j^2/2}$ converges and

$$E\left[\sum_{j=1}^{\infty} \frac{\lambda_j^4 \mu_j^2}{1 + \lambda_j^2}\right] = \sum_{j=1}^{\infty} \frac{\lambda_j^4 E(\mu_j^2)}{1 + \lambda_j^2} = \sum_{j=1}^{\infty} \lambda_j^2,$$

we have $E(U_1) = 1$ and $c_{d'} \to c$ (with $c$ a constant).

Arguing as in Proposition 4.1 in [3], we can show that

$$Var\left[\frac{1}{n} \sum_{i=1}^{n} U_i g(X_i)\right] \leq \frac{1}{n} E\left[U_1^2 g^2(X_1)\right] \to 0,$$

since a straightforward computation shows that $E(U_1^2) \leq M < \infty$ for all $d'$. Thus, under (i), the importance weights have the correct expectation and vanishing variance.

On the other hand, if (ii) holds, we can state the following proposition.

**Proposition 3.3.** *Under B, if $\sum_{j=1}^{\infty} \lambda_j^2 = \infty$ and $(\log n \log d')/\tau_{d'}^2 \to 0$, we have*

$$\frac{\tau_{d'}}{\sqrt{2 \log n}} E(T_{n,d'}) \to 1.$$

We note that our conditions imply that

$$\frac{\max_{1 \leq j \leq d'} \lambda_j^2 (1 + y_{0j}^2)}{\tau_{d'}^2} \to 0$$

so that asymptotic normality holds. The proof requires Lemmas A.1 and A.3.

The form reveals that it is possible to have much slower collapse than what Proposition 3.2 suggests. For instance, if $\lambda_j^2 = 1/j$, B holds but $\tau_{d'}^2 = \log d'(1 + o(1))$. In fact, the requirement that the $\lambda_j$ form an infinite sequence as above can be weakened to requiring simply that the $\lambda_j$ be bounded above uniformly, and this can be verified using a subsequence argument.

In conclusion, on the basis of Proposition 3.3, provided that the nonzero $\lambda_j$'s are commensurate, it seems reasonable to define $\sum_{j=1}^{d'} \lambda_j^2$ as the *effective dimension*. We note that the form of the effective dimension also plays a crucial role in the work of [7], who study Monte Carlo sample size requirements in the ensemble Kalman filter framework.

## Appendix

We first introduce two lemmas that pertain to Edgeworth expansion type uniform normal approximations of the distribution (the cdf and the density respectively) of independent sums of random variables. The two lemmas lay the groundwork for the proof of Proposition 3.1. Valid for moderately large deviations, the first result (Lemma A.1) is a special case of Theorem 2.5 in [12], and is stated here without proof.

**Lemma A.1.** *Let $\xi_1, \ldots, \xi_d$ be independent random variables with $E\xi_j = 0$ and $\sigma_j^2 = Var(\xi_j^2) < \infty$. Set*

$$S_d = \frac{1}{B_d}(\xi_1 + \cdots + \xi_d),$$

*where $B_d^2 = \sum_{j=1}^{d} \sigma_j^2$, and define the Lyapunov quotients*

$$L_{k,d} = \frac{1}{B_d^k} \sum_{j=1}^{d} E|\xi_j|^k, \quad k = 1, 2, \ldots.$$

*We also suppose $|E(Z_j^k)| \leq k! \gamma_j^{k-2} \sigma_j^2, k \geq 3$, where $\gamma_1, \ldots, \gamma_d$ are constant terms.*

*With these conditions, as $d \to \infty$, there exist analytic functions $P_d(x) = \sum_{k=3}^{\infty} \lambda_{k,d} x^k$ with $|\lambda_{k,d}| \leq Ac^k d^{-\frac{k-2}{2}}$ for some $A, c$ and all $d$, such that the cdf of $S_d$, denoted $G_d(\cdot)$, satisfies,*

$$1 - G_d(x) = (1 - \Phi(x)) exp(P_d(x))(1 + o(1)),$$

$$G_d(-x) = \Phi(-x) \exp(P_d(-x))(1 + o(1))$$

*uniformly for all $x \geq 0$ and $x = o(B_d/K_d)$, where $K_d = \max_{1 \leq j \leq d}\{\gamma_j, \sigma_j\}$. Furthermore, $P_d$ satisfies*

$$(9) \qquad |P_d(x)| \leq cx^3/B_d$$

*for some constant $c > 0$. We use $c$ generically as a constant independent of $d$.*

Lemma A.1 gives a normal approximation for the cdf of independent sums, and serves as the basis for the normality conditions of Proposition 3.1. Next we give a lemma for a normal approximation of the density of independent sums, which can be directly derived from Proposition 2 and Theorem 3 of [5].

**Lemma A.2.** *With the same notation and conditions as in Lemma A.1, we assume $\xi_{j,d}$ has density $g_{j,d}$ such that $\sup_x\{|g_{j,d}(x)| : 1 \leq j \leq d\} \leq M < \infty$. Then, as $d \to \infty$, the density of $S_d$, $g_d(\cdot) = G_d'(\cdot)$, satisfies*

$$g_d(x) = \phi(x)exp(P_d(x))(1 + o(1)),$$

$$g_d(-x) = \phi(-x)exp(P_d(-x))(1 + o(1))$$

*uniformly for all $x \geq 0$ and $x = o(B_d/K_d)$, where $K_d = \max_{1 \leq j \leq d}\{\gamma_j, \sigma_j\}$.*

We note in passing that the condition of uniform boundedness of the $g_{j,d}$ does not hold for $Z_j$, the Gaussian–Gaussian case. However, the sum of $\lambda_1^2 Z_1^2 + \lambda_2^2 Z_2^2$, where $\lambda_1, \lambda_2 > 0$ and $Z_1, Z_2$ are independent Gaussian, does indeed satisfy the condition. This may be verified by a direct calculation of the density of the convolution.

The next lemma is given for the purpose of verifying the Lyapunov quotients conditions appearing in Lemmas A.1 and A.2.

**Lemma A.3.** *Let $Z_j, V_j, \epsilon_j, j = 1, \ldots, d$, be iid $N(0,1)$. Let $\lambda_1 \geq \lambda_2 \geq \cdots$ where $\sum_{j=1}^{\infty} \lambda_j^2 = \infty$. Then, given $\mu_j \equiv V_j + \frac{\epsilon_j}{\lambda_j}$, for all $j$, we have*

$$(10) \qquad \begin{aligned} &\lambda_j^{2k} E\big(|(Z_j + \mu_j)^2 - (1 + \mu_j^2)|^k \big| \mu_j\big) \\ &\leq \frac{O_p(\sqrt{\log d})^k}{k!} \rho^k \lambda_j^4 E\big((Z_j + \mu_j)^2 - (1 + \mu_j^2)\big| \mu_j\big), \end{aligned}$$

*for $k \geq 3$.*

Thus, given the mean vector $\xi = [\mu_1, \mu_2, \cdots, \mu_d]$ defined in (4), Lemma A.3 states that the Lyapunov conditions required by Lemma A.1 hold, with probability tending to 1. We note that our argument also implies Lemma A.1.

*Proof of Lemma A.3.* Since $(Z_j + \mu_j)^2 - (1 + \mu_j^2) = (Z_j^2 - 1) + 2\mu_j Z_j$, it is enough to bound

$$\lambda_j^{2k} E\big(|(Z_j^2 - 1) + 2\mu_j Z_j|^k \big| \mu_j\big) \leq 2^k\big(\lambda_j^{2k} E|Z_j^2 - 1|^k + 2^k(|\mu_j|\lambda_j^2)^k E|Z_j|^k\big).$$

By standard properties of the Gaussian moments, for some positive constant $C$,

$$E|Z_j^2 - 1|^k \leq C^k k!, \quad \text{and} \quad E|Z_j|^k \leq C^k k! \ .$$

Since $E(Z_j^2 - 1 + 2\mu_j Z_j)^2 = 2 + 4\mu_j^2$ we see that (10) follows from the bound

$$\lambda_j^{2k}|\mu_j|^k \leq \lambda_j^2 \mu_j^2 \max\{|\lambda_\ell^2 \mu_\ell|^{k-2} : 1 \leq \ell \leq d\} = \big(O_p(\sqrt{\log d})\big)^{k-2},$$

since the $\lambda_j^2 \mu_j$ are independent $N(0, \lambda_j^2 + \lambda_j^4)$ so that

$$\max\{|\lambda_\ell^2 \mu_\ell| : 1 \leq \ell \leq d\} \leq (\lambda_j^2 + \lambda_j^4)^{1/2} \max\{|V_\ell| : 1 \leq \ell \leq d\}$$

where the $V_\ell$ are i.i.d. $N(0,1)$. The lemma follows. $\qquad \square$

The remainder of the Appendix is devoted to the proof of the main result given in Proposition 3.1.

*Proof of Proposition 3.1.* Let $S_j$ $(j = 1, \dots, n)$ be as defined in the Proposition and let $S_{(1)}$ be the minimum. Note that

$$(11) \qquad E(T_{n,d}|S_{(1)}) = \frac{(n-1)\int_{S_{(1)}}^{\infty} \exp\big(-\tau_d(z - S_{(1)})\big)\mathrm{d}G_d(z)}{\bar{G}_d(S_{(1)})},$$

since, given $S_{(1)}$, the remaining $(n-1)$ observations are i.i.d. with cdf equal to $G_d(z)/\bar{G}_d(S_{(1)})$, $z \geq S_{(1)}$.

Let $\varepsilon_d$ be a sequence of constants such that $\varepsilon_d \to 0$ and $\varepsilon_d\tau_d/\sqrt{2\log n} \to \infty$ as $n, d \to \infty$. We first define, for $x < \varepsilon_d\tau_d$,

$$(12) \qquad h_{n,d}(x) := \int_x^{\infty} \exp\big(-\tau_d(z - x)\big)\mathrm{d}G_d(z).$$

To evaluate $h_{n,d}(x)$, we break the integral into two parts: the first part yields the integral from $x$ to $x + \varepsilon_d\tau_d$, and the second part yields the tail integral from $x + \varepsilon_d\tau_d$ to $\infty$. By using the normal approximations of Lemmas A.1 and A.2, under the assumption that $(\log n)/\tau_d^2 \to 0$, one can show that the second part is $o\big(\sqrt{2\log n}/n\tau_d\big)$.

To deal with the first part, we shall show that as $x \to -\infty$ and $x > -\varepsilon_d\tau_d$,

$$(13) \qquad \int_x^{x+\varepsilon_d\tau_d} \exp\big(-\tau_d(z - x)\big)\mathrm{d}G_d(z) = \frac{1}{\tau_d}\phi(x)\exp\big(P_d(x)\big)(1 + o(1))$$

To this end, applying Lemma A.2 with $\ell = 3$, we obtain,

$$
\begin{aligned}
R_d(x) &:= \int_x^{x+\varepsilon_d\tau_d} \exp\Big[-\tau_d(z - x) - \frac{1}{2}(z^2 - x^2) + P_d(z) - P_d(x)\Big]\mathrm{d}z(1 + o(1)) \\
&= \int_0^{\Delta_{n,d}} \exp\Big[-\tau_d v - \frac{1}{2}((x+v)^2 - x^2) \\
&\qquad\qquad + \sum_{k=3}^{\infty} \lambda_{k,d}((x+v)^k - x^k)\Big]\mathrm{d}v(1 + o(1)) \\
&= \frac{1}{|x|}\int_0^{|x|\varepsilon_d\tau_d} \exp\Big[-(-1 + \frac{\tau_d}{|x|})w - \frac{w^2}{2|x|^2} \\
&\qquad\qquad + \sum_{k=3}^{\infty} \lambda_{k,d}\sum_{j=1}^{k}(-1)^{k-j}C_{k,j}|x|^{k-2j}w^j\Big]\mathrm{d}w(1 + o(1)) \\
(14) \quad &= \frac{1}{|x|}\int_0^{|x|\varepsilon_d\tau_d} \exp\Big[-b_1 w - b_2 w + \sum_{j=3}^{\infty} b_j w^j\Big]\mathrm{d}w(1 + o(1)),
\end{aligned}
$$

where

$$
\begin{aligned}
b_1 &= \frac{\sigma\sqrt{d}}{|x|} - b_1^* = -1 + \frac{\tau_d}{|x|} - \sum_{k=3}^{\infty}(-1)^{k-1}C_{k,1}\lambda_{k,d}|x|^{k-2}, \\
b_2 &= \frac{1}{2|x|^2} - b_2^* = \frac{1}{2|x|^2} - \sum_{k=3}^{\infty}(-1)^{k-2}C_{k,2}\lambda_{k,d}|x|^{k-4}
\end{aligned}
$$

and

$$b_j = \sum_{k=j}^{\infty} (-1)^{k-j} C_{k,j} \lambda_{k,d} |x|^{k-2j}.$$

Note $|\lambda_{k,d}| \leq A c_0^k \tau_d^{-(k-2)}$ and $C_{k,j} < c^k$, for some constants $A, c_0, c$ where $\mu_j = V_j + \epsilon_j/\lambda_j$. Hereafter, we use $c$ as a generic positive constant that does not depend on $x$ and $d$. Under the assumptions that $x \to -\infty$, $|x| < \varepsilon_d \tau_d$ (hence $|x|/\tau_d \to 0$), and $|x|\Delta_{n,d} \to \infty$, we have, firstly,

$$
\begin{aligned}
b_1^* &\leq \sum_{k=3}^{\infty} A k c_0^k \tau_d^{-(k-2)} |x|^{k-2} \\
&= A c_0^2 [3(c_0|x|/\tau_d)/(1 - (c_0|x|/\tau_d)) + (c_0|x|/\tau_d)^2/(1 - (c_0|x|/\tau_d))^2] \\
&= o(1),
\end{aligned}
\tag{15}
$$

secondly,

$$b_2^* \leq x^{-2} \sum_{k=3}^{\infty} c(c|x|/\tau_d)^{k-2} = |x|^{-2}(c|x|/\tau_d)/(1 - c|x|/\tau_d) = o(|x|^{-2}), \tag{16}$$

and thirdly,

$$
\begin{aligned}
b_j &\leq (\sum_{k=j}^{2j-1} + \sum_{k=2j}^{\infty}) A(cc_0)^k \tau_d^{-(k-2)} |x|^{k-2j} \\
&= \sum_{k=j}^{2j-1} A(cc_0)^k (|x|/\tau_d)^{k-j} |x|^{-j} \tau_d^{-(j-2)} \\
&\quad + \sum_{k=2j}^{\infty} A(cc_0)^k (|x|/\tau_d)^{k-2j} \tau_d^{2-2j} \\
&\leq |x|^{-2}(c|x|\tau_d)^{-(j-2)} + c\tau_d^{2-2j} \\
&\leq 2|x|^{-2}(c|x|\tau_d)^{-(j-2)}.
\end{aligned}
\tag{17}
$$

Since $w/(|x|\tau_d) \leq \varepsilon_d \to 0$, we can further derive

$$
\begin{aligned}
\sum_{j=3}^{\infty} b_j w^j &\leq 2(w/x)^2 (cw/(|x|\tau_d))^{j-2} = 2(w/|x|)^2 cw/(|x|\tau_d)/[1 - cw/(|x|\tau_d)] \\
&= o(|x|^{-2})w^2.
\end{aligned}
$$

Combining (14), (15), (16), and (18) yields

$$R_d(x) = \frac{1}{|x|} \int_0^{|x|\Delta_{n,d}} \exp\big[-(-1 + \frac{\tau_d}{|x|})(1 + o(1))w - (\frac{w^2}{2|x|^2})(1 + o(1))\big] \mathrm{d}w. \tag{18}$$

The $o(1)$'s appearing in the last expression are uniform as $w$ varies over the integral interval. Now, the bounded convergence theorem ensures $R_d(x) = (1/\tau_d)(1 + o(1))$, which establishes (13). Taking into account the remainder term, we conclude that

$$h_{n,d}(x) = \frac{1}{\tau_d} \phi(x) \exp\big(P_d(x)\big)(1 + o(1)) + o\big(\frac{\sqrt{2\log n}}{n\tau_d}\big). \tag{19}$$

Our target $(\tau_d/\sqrt{2\log n})E(T_{n,d})$ can now be written as

$$
\begin{aligned}
\frac{\tau_d}{\sqrt{2\log n}}E(T_{n,d}) &= \frac{\tau_d(n-1)}{\sqrt{2\log n}}E\big[\frac{h_{n,d}(S_{(1)})}{\bar{G}_d(S_{(1)})}\big] \\
&= \frac{\tau_d n}{\sqrt{2\log n}}\int_{-\infty}^{\infty}h_{n,d}(x)\bar{G}_d^{n-2}(x)dG_d(x).
\end{aligned}
\tag{20}
$$

We decompose the preceding integral into three parts

$$
\frac{\tau_d}{\sqrt{2\log n}}E(T_{n,d}) = I_{n,d} + II_{n,d} + III_{n,d}
\tag{21}
$$

where $I_{n,d}, II_{n,d},$ and $III_{n,d}$ represent the integral of (11) over the intervals $[-\infty, -\varepsilon_d\tau_d], (-\varepsilon_d\tau_d, -(\log n)^{1/4}),$ and $[-(\log n)^{1/4}, \infty),$ respectively. The preceding discussion, combined with the approximation $g_d(x) = xG_d(x)(1+o(1))$ as $x \to -\infty$ and $|x| = o(\tau_d)$, implies that the dominating part is the quantity represented by $II_{n,d}$. We have,

$$
\begin{aligned}
II_{n,d} &= \frac{n(n-1)}{\sqrt{2\log n}}\int_{-\varepsilon_d\tau_d}^{-(\log n)^{1/4}}xG_d(x)\bar{G}_d^{n-2}(x)dG_d(x)(1+o(1)) \\
&= \frac{1}{\sqrt{2\log n}}\int_{nG_d(-\varepsilon_d\tau_d)}^{nG_d(-(\log n)^{1/4})}G_d^{-1}(w/n)w(1-w/n)^n dw(1+o(1)) \\
&= \frac{1}{\sqrt{2\log n}}\int_{nG_d(-\varepsilon_d\tau_d)}^{nG_d(-(\log n)^{1/4})}\sqrt{-2\log(w/n)}we^{-w}dw(1+o(1)) \\
&= \int_0^{\infty}we^{-w}dw(1+o(1)) - \frac{1}{\sqrt{2\log n}}\int_0^{\infty}w\log we^{-w}dw(1+o(1)) \\
&= 1+o(1).
\end{aligned}
\tag{22}
$$

To arrive at (22) we have used the approximation $G_d^{-1}(z) = \sqrt{-2\log z}(1+o(1))$ for $z \to 0$ in light of Lemma A.1 and Mill's ratio.

For the remaining two parts, we use Mill's ratio and obtain

$$
\begin{aligned}
I_{n,d} + III_{n,d} &\leq \frac{\tau_d}{\sqrt{2\log n}}(n-1)\big[P(S_{(1)} \leq -\varepsilon_d\tau_d) + P(S_{(1)} \geq -(\log n)^{1/4})\big] \\
&= \frac{\tau_d}{\sqrt{2\log n}}(n-1)\big[1 - \bar{G}_d^n(-\varepsilon_d\tau_d) + \bar{G}_d^n\big(-(\log n)^{1/4}\big)\big] \\
&\leq \frac{\tau_d}{\sqrt{2\log n}}(n-1)\big[nG_d(-\varepsilon_d\tau_d) + \bar{G}_d^n\big(-(\log n)^{1/4}\big)\big] \\
&\to 0.
\end{aligned}
\tag{23}
$$

Finally, combining (21), (22), and (23), yields the desired result. $\qquad\square$